\documentclass[journal,10pt]{IEEEtran}
\usepackage[named]{algo}
\usepackage{graphics}
\usepackage{graphicx}
\usepackage{stfloats}
\usepackage{amsmath}
\usepackage{amsfonts}
\usepackage{subfigure}
\usepackage{float}
\usepackage{url}

\newcommand{\Section}{\section}
\newcommand{\SubSection}{\subsection}

\newtheorem{theorem}{Theorem}

\newcommand{\bbeta}{\mbox{\boldmath  $\beta$}}

\newcommand{\bPsi}{\mbox{\boldmath   $\Psi$}}

\newcommand{\blambda}{\mbox{\boldmath $\lambda$}}

\def\urltilda{\kern -.15em\lower .7ex\hbox{\~{}}\kern .04em}
\def\urldot{\kern -.10em.\kern -.10em}
\def\urlhttp{http\kern -.10em\lower -.1ex\hbox{:}\kern -.12em\lower 0ex\hbox{/}\kern -.18em\lower 0ex\hbox{/}}

\begin{document}

\title{Fast Image Recovery Using Variable
Splitting \\ and Constrained Optimization}

\author{Manya V. Afonso, Jos\'{e} M. Bioucas-Dias, and M\'{a}rio A. T. Figueiredo
\thanks{The authors are with the {\it Instituto de Telecomunica\c{c}\~{o}es} and the
Department of Electrical and Computer Engineering,
{\it Instituto Superior T\'{e}cnico,} 1049-001 Lisboa, {\bf Portugal}.}
\thanks{Emails: mafonso@lx.it.pt, bioucas@lx.it.pt, mario.figueiredo@lx.it.pt}
\thanks{A preliminary much shorter version of this work appeared in \cite{FigueiredoDiasAfonsoSSP09}.}}

\maketitle

\begin{abstract}
We propose a new fast algorithm for solving one of the standard
formulations of image restoration and reconstruction which consists of
an unconstrained optimization problem where the objective
includes an $\ell_2$ data-fidelity term and a non-smooth
regularizer. This formulation allows both wavelet-based (with
orthogonal or frame-based representations) regularization or
total-variation regularization. Our approach is based on a
variable splitting to obtain an equivalent constrained
optimization formulation, which is then addressed with an
augmented Lagrangian method. The proposed algorithm is an
instance of the so-called {\it alternating direction method
of multipliers}, for which convergence has been proved.
Experiments on a set of image restoration and reconstruction
benchmark problems show that the proposed algorithm is faster
than the current state of the art methods.

\end{abstract}

\Section{Introduction}
\SubSection{Problem Formulation}
\label{sec:intro}
Image restoration/reconstruction is one of the earliest and most classical linear
inverse problems in imaging, dating back to the 1960's
\cite{Andrews}. In this class of problems,  a noisy indirect observation ${\bf y}$,
of an original image ${\bf x}$, is modeled as
\[
{\bf y} = {\bf B}{\bf x} + {\bf n},
\]
where ${\bf B}$  is the matrix representation of the direct
operator and ${\bf n}$ is noise. As is common, we are adopting
the vector notation for images, where the pixels on an
$M\times N$ image are stacked into a an $(NM)$-vector in,
e.g., lexicographic order. In the sequel, we denote by $n$
the number of elements of ${\bf x}$, thus ${\bf x}\in\mathbb{R}^n$,
while  ${\bf y}\in\mathbb{R}^m$ ($m$ and $n$ may or may not be equal).

In the particular case of image deblurring/deconvolution, ${\bf B}$
is the matrix representation of a convolution operator; if this
convolution is periodic, ${\bf B}$ is then a (block) circulant
matrix. This type of observation model describes well several
physical mechanisms, such as relative motion between the camera and
the subject (motion blur), bad focusing (defocusing blur), or a
number of other mechanisms which are well modeled by a convolution.

In more general image reconstruction problems, ${\bf B}$ represents
some linear direct operator, such as a set of tomographic projections
(Radon transform), a partially observed (e.g., Fourier) transform, or the loss
of part of the image pixels.

It is well known that the problem of estimating ${\bf x}$ from ${\bf y}$
is ill-posed, thus this inverse problem can only be solved satisfactorily
by adopting some sort of regularization (or prior information, in Bayesian
inference terms).  One of the standard formulations for wavelet-based
regularization of image restoration/reconstruction problems is built as
follows. Let the unknown image ${\bf x}$ be represented as a linear
combination of the elements of some frame, {\it i.e.},
${\bf x} = {\bf W}\bbeta$, where $\bbeta\in \mathbb{R}^d$, and the
columns of the $n\times d$ matrix ${\bf W}$ are the elements of a
wavelet\footnote{We will use the generic term ``wavelet" to mean any wavelet-like
multi-scale representation, such as ``curvelets", ``beamlets", or ``ridgelets".}
frame (an orthogonal basis or a redundant dictionary).
Then, the coefficients of this representation are estimated
from the noisy image, under one of the well-known
sparsity inducing regularizers, such as the $\ell_1$ norm
(see \cite{DaubechiesDefriseDeMol}, \cite{EladCVPR2006},
\cite{FigueiredoDiasNowak2007}, \cite{FigueiredoNowak2003},
\cite{FigueiredoNowak2005}, and further references therein).
Formally, this leads to the following optimization problem:
\begin{equation}
\widehat{\bbeta} = \arg\min_{\bbeta} \frac{1}{2}\| {\bf BW}\bbeta -
{\bf y}\|_2^2 + \tau \, \phi (\bbeta) \label{problem}
\end{equation}
where $\phi:\mathbb{R}^d \! \rightarrow \bar{\mathbb{R}}$, usually
called the {\it regularizer} or {\it regularization function} is
usually nonsmooth, or maybe even nonconvex, and $\tau \geq 0$ is
the regularization parameter. This formulation is referred to as
the {\it synthesis approach} \cite{EladMilanfarRubinstein}, since
it is based on a synthesis equation where ${\bf x}$ is
synthesized from its representation coefficients
(${\bf x} = {\bf W}\bbeta$) with are the object of the
estimation criterion. Of course, the final image estimate is
computed as $\widehat{\bf x} = {\bf W}\widehat{\bbeta} $.

An alternative formulation applies a regularizer directly
to the unknown image, leading to criteria of the form
\begin{equation}
\widehat{{\bf x}} = \arg\min_{{\bf x}} \frac{1}{2}\| {\bf B\, x} -
{\bf y}\|_2^2 + \tau \, \phi ({\bf x}) \label{problem2}
\end{equation}
where $\phi:\mathbb{R}^n \! \rightarrow \bar{\mathbb{R}}$ is the
regularizer. This type of criteria are usually called {\it analysis
approaches}, since they're based on a regularizer that
analyzes the image itself, $\phi({\bf x})$, rather than
the coefficients of a representation thereof.
Arguably, the best known and most often used
regularizer used in analysis approaches to image restoration is
the total variation (TV) norm \cite{ROF},  \cite{Chan_et_al_2005}.
Wavelet-based analysis approaches are also possible \cite{EladMilanfarRubinstein},
but will not be considered in this paper.


Finally, it should be mentioned that problems (\ref{problem})
and (\ref{problem2}) can be seen as the Lagrangians
of associated constrained optimization problems: (\ref{problem}) is
the Lagrangian of the constrained problem
\begin{equation}
\label{eq:bp}
\min_{\bbeta} \phi(\bbeta)\hspace{0.5cm} \mbox{subject to}
\hspace{0.5cm} \| {\bf BW}\bbeta - {\bf y}\|_2^2 \leq
\varepsilon,
\end{equation}
while (\ref{problem2}) is the Lagrangian
of
\begin{equation}
\label{eq:bp2}
\min_{{\bf x}} \phi({\bf x})\hspace{0.5cm} \mbox{subject to}
\hspace{0.5cm} \| {\bf B\, x}- {\bf y}\|_2^2 \leq
\varepsilon.
\end{equation}
Specifically, a solution of (\ref{eq:bp}) (for any $\varepsilon$
such that this problem is feasible) is either the null vector,
or else is a minimizer of (\ref{problem}), for some $\tau>0$
(see \cite[Theorem~27.4]{Roc70}). A similar relationship
exists between problems (\ref{problem2}) and (\ref{eq:bp2}).

\SubSection{Previous Algorithms}
For any problem of non-trivial
dimension, matrices ${\bf BW}$, ${\bf B}$, and ${\bf W}$
cannot be stored explicitly, and it is costly, even impractical,
to access portions (lines, columns, blocks) of them.
On the other hand,
matrix-vector products involving ${\bf B}$ or ${\bf W}$
(or their conjugate transposes ${\bf B}^H$ and ${\bf W}^H$) can
be done quite efficiently. For example, if the columns of
${\bf W}$ contain a wavelet basis or a tight wavelet frame,
any multiplication of the form ${\bf W}{\bf v}$ or ${\bf W}^H{\bf v}$
can be performed by a fast wavelet transform algorithm \cite{Mallat}.
Similarly, if ${\bf B}$ represents a convolution, products of the form
${\bf B}{\bf v}$ or ${\bf B}^H{\bf v}$ can be performed with the help
of the fast Fourier transform (FFT) algorithm. These facts have
stimulated the development of special purpose
methods, in which the only operations involving ${\bf B}$ or
${\bf W}$ (or their conjugate transposes) are matrix-vector products.

To present a unified view of algorithms for handling (\ref{problem}) and
(\ref{problem2}), we write them in a common form
\begin{equation}
\min_{{\bf x}} \frac{1}{2}\| {\bf A\, x} -
{\bf y}\|_2^2 + \tau \, \phi ({\bf x}) \label{problem_common}
\end{equation}
where ${\bf A}={\bf BW}$, in the case of (\ref{problem}), while
${\bf A}={\bf B}$, for (\ref{problem2}).

Arguably, the standard algorithm for solving problems of the form
(\ref{problem_common}) is the so-called {\it iterative shrinkage/thresholding} (IST)
algorithm. IST can be derived as an expectation-maximization
(EM) algorithm \cite{FigueiredoNowak2003}, as a {\it majorization-minimization} (MM, \cite{Hunter})
method \cite{DaubechiesDefriseDeMol}, \cite{FigueiredoNowak2005},
or as a forward-backward splitting technique \cite{CombettesSIAM}, \cite{Hale}.
A key ingredient of IST algorithms is the so-called shrinkage/thresholding
function, also known as the Moreau proximal mapping \cite{CombettesSIAM} or the denoising function,
 associated to the
regularizer $\phi$, which provides the solution
of the corresponding pure denoising problem. Formally, this function is denoted
as $\bPsi_{\tau\phi}:\mathbb{R}^m\rightarrow \mathbb{R}^m$ and
defined as
\begin{equation}
\bPsi_{\tau\phi}({\bf y}) = \arg\min_{{\bf x}} \frac{1}{2}\|{\bf x}-{\bf y}\|_2^2 + \tau\phi({\bf x}).
\label{MPM}
\end{equation}
Notice that if $\phi$ is proper and convex, the function being minimized is
proper and strictly convex, thus the minimizer exists and is unique making the
function well defined \cite{CombettesSIAM}.

For some choices of $\phi$, the corresponding denoising functions $\bPsi_{\tau\phi}$
have well known closed forms. For example, choosing $\phi({\bf x}) = \|{\bf x}\|_1 = \sum_i |x_i|$,
the $\ell_1$ norm, leads to $\bPsi_{\tau\ell_1}({\bf y}) = \mbox{soft}({\bf y}, \tau)$, where
$\mbox{soft}(\cdot, \tau)$ denotes the component-wise application of the function $y \mapsto \mbox{sign}(y)\max\{|y|-\tau,0\}$.

If $\phi({\bf x}) = \|{\bf x}\|_0 = |\{i: x_i \neq 0\}|$, usually referred to as the
$\ell_0$ ``norm" (although it is not a norm), despite the fact that this regularizer
is not convex, the corresponding shrinkage/thresholding function also has a simple
close form: the so-called hard-threshold function, $\bPsi_{\tau\ell_0}({\bf y}) =
\mbox{hard}({\bf y}, \sqrt{2\, \tau})$, where  $\mbox{hard}(\cdot, a)$ denotes the
component-wise application of the function $y \mapsto y 1_{|y|\geq a}$.
A comprehensive coverage of {\it Moreau proximal maps} can be found in \cite{CombettesSIAM}.

Each IST iteration for solving (\ref{problem_common}) is given by
\begin{equation}
{\bf x}_{k+1} = \bPsi_{\tau\phi}\left({\bf x}_t - \frac{1}{\gamma}\, {\bf A}^H\left(
{\bf A} {\bf x}_k - {\bf y}\right)\right),
\end{equation}
where $1/\gamma$ is a step size. Notice that ${\bf A}^H\left({\bf A} {\bf x}_k - {\bf y}\right)$
is the gradient of the data-fidelity term $(1/2)\|{\bf A} {\bf x} - {\bf y}\|_2^2$,
computed at ${\bf x}_k$; thus, each IST iteration takes a step of length $1/\gamma$ in
the direction of the negative gradient of the data-fidelity term, followed by the application of the shrinkage/thresholding function associated with the regularizer $\phi$.

It has been shown that if $\gamma > \|{\bf A}\|_2^2/2$
and $\phi$ is convex, the algorithm converges to a solution of (\ref{problem})
\cite{CombettesSIAM}. However, it is known that IST may be quite slow, specially
when $\tau$ is very small and/or the matrix ${\bf A}$ is very ill-conditioned
\cite{FISTA}, \cite{TwIST}, \cite{FigueiredoDiasNowak2007}, \cite{Hale}.
This observation has stimulated
work on faster variants of IST, which we will briefly review in the next
paragraphs.

In the {\it two-step IST} (TwIST) algorithm \cite{TwIST}, each iterate depends
on the two previous iterates, rather than only on the previous one
(as in IST). This algorithm may be seen as a non-linear version
of the so-called two-step methods for linear problems \cite{Axelsson}.
TwIST was shown to be considerably faster than IST on a
variety of wavelet-based and TV-based image restoration problems;
the speed gains can reach up two orders of magnitude in typical benchmark
problems.

Another two-step variant of IST, named {\it fast IST algorithm} (FISTA),
was recently proposed and also shown to clearly outperform IST in terms of
speed \cite{FISTA}. FISTA is a non-smooth variant of Nesterov's
optimal gradient-based algorithm for smooth convex problems \cite{Nesterov},
\cite{Nesterov_Book}.

A strategy recently proposed to obtain faster variants of IST consists in relaxing the condition
$\gamma > \gamma_{\mbox{\tiny min}} \equiv \|{\bf A}\|_2^2/2$. In the
SpaRSA (standing for {\it sparse reconstruction by separable approximation})
framework \cite{SpaRSA_ICASSP}, \cite{SpaRSA_SP},
a different $\gamma_t$ is used in each iteration (which may be smaller
than  $\gamma_{\mbox{\tiny min}}$, meaning larger step sizes). It
was shown experimentally that SpaRSA clearly outperforms standard IST.
A convergence result for SpaRSA was also given in \cite{SpaRSA_SP}.

Finally, when the slowness is caused by the use of a small value of the
regularization parameter, {\it continuation} schemes have been found quite
effective in speeding up the algorithm. The key observation is that IST
algorithm benefits significantly from {\it warm-starting}, {\it i.e.}, from being
initialized near a minimum of  the objective function. This suggests
that we can use the solution of (\ref{problem_common}), for a given value of
$\tau$, to initialize IST in solving the same problem for a nearby
value of $\tau$.  This {\it warm-starting} property underlies  {\it continuation}
schemes \cite{GPSR}, \cite{Hale}, \cite{SpaRSA_SP}. The idea is to use
IST to solve (\ref{problem})  for a larger value of $\tau$ (which is usually fast),
then decrease $\tau$ in steps toward its desired value, running IST with warm-start
for each successive value of $\tau$.

\SubSection{Proposed Approach}
The approach proposed in this paper is based on the principle of variable
splitting, which goes back at least to Courant in the 40's \cite{Courant},
\cite{Wang}. Since the objective function (\ref{problem_common}) to be minimized
is the sum of two functions, the idea is to split the variable ${\bf x}$
into a pair of variables, say ${\bf x}$ and ${\bf v}$, each to serve as the
argument of each of the two functions, and then minimize the sum of
the two functions under the constraint that the two variables have to be equal,
so that the problems are equivalent. Although variable splitting is also the
rationale behind the recently proposed split-Bregman method \cite{GoldsteinOsher},
in this paper, we exploit a different type of splitting to attack
problem (\ref{problem_common}). Below we will explain this difference in detail.

The obtained constrained optimization problem is then dealt with using an augmented
Lagrangian (AL) scheme \cite{NocedalWright}, which is known to be equivalent
to the Bregman iterative methods recently proposed to handle imaging inverse
problems (see \cite{YinOsherGoldfarbDarbon} and references therein). We prefer the
AL perspective, rather than the Bregman iterative view, as it is a standard
and more elementary optimization tool (covered in most textbooks on optimization).
In particular, we solve the constrained problem resulting from the variable splitting
using an algorithm known as alternating direction method of multipliers (ADMM)
\cite{EcksteinBertsekas}.

The application of ADMM to our particular problem involves solving a linear system
with the size of the unknown image (in the case of problem (\ref{problem2})) or with
the size of its representation (in the case of problem (\ref{problem})). Although this
seems like an unsurmountable obstacle, we show that it is not the case. In many problems
of the form (\ref{problem2}), such as deconvolution, recovery of missing samples, or
reconstruction from partial Fourier observations, this system can be solved very quickly
in closed form (with $O(n)$ or $O(n\log n)$ cost). For problems of the form (\ref{problem}),
we show how exploiting the fact that ${\bf W}$ is a tight Parseval frame, this system can still
be solved efficiently (typically with $O(n\log n)$ cost.

We report results of a comprehensive set of experiments, on a set of benchmark problems,
including image deconvolution, recovery of missing pixels, and reconstruction from partial
Fourier transform, using both frame-based and TV-based regularization. In all the experiments,
the resulting algorithm is consistently and considerably faster than the previous state of the
art methods FISTA \cite{FISTA}, TwIST \cite{TwIST}, and SpaRSA \cite{SpaRSA_SP}.

The speed of the proposed algorithm, which we term SALSA ({\it split augmented
Lagrangian shrinkage algorithm}), comes from the fact that it uses (a regularized version of)
the Hessian of the data fidelity term of (\ref{problem_common}), that is, ${\bf A}^H{\bf A}$,
while the above mentioned algorithms essentially only use gradient information.

\subsection{Organization of the Paper}
Section \ref{sec:tools} describes the basic ingredients of SALSA: variable splitting, augmented Lagrangians,
and ADMM. In Section \ref{sec:salsa}, we show how these ingredients are combined to obtain
the proposed SALSA. Section \ref{sec:experiments} reports experimental results, and
Section \ref{sec:conclusions} ends the paper with a few remarks and pointers to future work.

\Section{Basic Ingredients}

\label{sec:tools}
\SubSection{Variable Splitting}
Consider an unconstrained optimization problem in which
the objective function is the sum of two functions, one
of which is written as the composition of two functions,
\begin{equation}
\min_{{\bf u}\in \mathbb{R}^n} f_1({\bf u}) + f_2\left(g({\bf u})\right),\label{unconstrained_basic}
\end{equation}
where $g:\mathbb{R}^n\rightarrow \mathbb{R}^d$.
Variable splitting is a very simple procedure that consists in
creating a new variable, say ${\bf v}$,
to serve as the argument of $f_2$, under the
constraint that $g({\bf u}) = {\bf v}$. This leads to the constrained problem
\begin{equation}\begin{array}{cl}
 {\displaystyle \min_{{\bf u}\in \mathbb{R}^n,\, {\bf v}\in\mathbb{R}^d}} & f_1({\bf u}) + f_2({\bf v})\\
 \mbox{subject to} &g({\bf u}) = {\bf v}, \end{array}\label{constrained_basic}
\end{equation}
which is clearly equivalent to unconstrained problem (\ref{unconstrained_basic}):
in the feasible set $\{({\bf u},{\bf v}): g({\bf u}) = {\bf v}\}$, the objective
function in (\ref{constrained_basic}) coincides with that in (\ref{unconstrained_basic}).
The rationale behind variable splitting methods is that it may
be easier to solve the constrained problem (\ref{constrained_basic})
than it is to solve its unconstrained counterpart (\ref{unconstrained_basic}).

The splitting idea has been recently used in several image processing applications.
A variable splitting method was used in \cite{Wang} to obtain a fast
algorithm for TV-based image restoration. Variable splitting was also
used in \cite{BioucasFigueiredo2008} to handle problems involving compound
regularizers; {\it i.e.}, where instead of a single regularizer
$\tau \phi({\bf x})$ in (\ref{problem_common}), one has a linear combination of
two (or more) regularizers $\tau_1 \phi_1({\bf x}) + \tau_2 \phi_2({\bf x})$.
In \cite{BioucasFigueiredo2008} and \cite{Wang}, the constrained problem
(\ref{constrained_basic}) is attacked by a quadratic penalty approach, i.e.,
by solving
\begin{equation}
 {\displaystyle \min_{{\bf u}\in \mathbb{R}^n,\, {\bf v}\in\mathbb{R}^d}}\;\;
   f_1({\bf u}) + f_2({\bf v}) +
 \frac{\alpha}{2}\, \|g({\bf u}) -{\bf v}\|_2^2,
\label{quadratic_penalty}
\end{equation}
by alternating minimization with respect to ${\bf u}$ and ${\bf v}$, while
slowly taking $\alpha$ to very large values (a {\it continuation} process),
to force the solution of (\ref{quadratic_penalty}) to approach that
of (\ref{constrained_basic}), which in turn is equivalent to (\ref{unconstrained_basic}).
The rationale behind these methods is that each step of this
alternating minimization may be much easier than the original
unconstrained problem (\ref{unconstrained_basic}). The drawback is that
as $\alpha$ becomes very large, the intermediate minimization problems
become increasingly ill-conditioned, thus causing numerical problems
(see \cite{NocedalWright}, Chapter 17).

A similar variable splitting approach underlies the recently
proposed split-Bregman methods \cite{GoldsteinOsher}; however,
instead of using a quadratic penalty technique, those methods
attack the constrained problem directly using a Bregman iterative
algorithm \cite{YinOsherGoldfarbDarbon}. It has been shown that,
when $g$ is a linear function, {\it i.e.}, $g({\bf u}) = {\bf Gu}$,
the Bregman iterative algorithm is equivalent to the augmented
Lagrangian method \cite{YinOsherGoldfarbDarbon}, which is briefly
reviewed in the following subsection.

\SubSection{Augmented Lagrangian}
Consider the constrained optimization problem
\begin{equation}\begin{array}{cl}
 {\displaystyle \min_{{\bf z}\in \mathbb{R}^n}} & E({\bf z})\\
 \mbox{s.t.} & {\bf A z - b = }\mbox{\boldmath $0$}, \end{array}\label{constrained_linear}
\end{equation}
where ${\bf b} \in \mathbb{R}^p$ and ${\bf A}\in \mathbb{R}^{p\times n}$,
{\it i.e.}, there are $p$ linear equality constraints.
The so-called augmented Lagrangian function for this problem is
defined as
\begin{equation}
{\cal L}_A ({\bf z},\blambda,\mu) = E({\bf z}) + \blambda^T ({\bf b-Az}) +
\frac{\mu}{2}\,  \|{\bf Az-b}\|_2^2,\label{augmented_L}
\end{equation}
where $\blambda \in \mathbb{R}^p$ is a vector of Lagrange multipliers
and $\mu \geq 0$ is called the penalty parameter \cite{NocedalWright}.

The so-called {\it augmented Lagrangian method} (ALM) \cite{NocedalWright},
also known as the {\it method of multipliers} (MM) \cite{Hestenes}, \cite{Powell},
consists in minimizing ${\cal L}_A ({\bf z},\blambda,\mu)$
with respect to ${\bf z}$, keeping $\blambda$ fixed, then updating
$\blambda$, and repeating these two steps until some convergence criterion
is satisfied. Formally, the ALM/MM works as follows:

\vspace{0.3cm}
\begin{algorithm}{ALM/MM}{
\label{alg:salsa1}}
Set $k=0$, choose $\mu > 0$, ${\bf z}_0$, and  $\blambda_0$.\\
\qrepeat\\
     ${\bf z}_{k+1} \in \arg\min_{{\bf z}} {\cal L}_A ({\bf z},\blambda_k,\mu)$\\
     $\blambda_{k+1} = \blambda_{k} + \mu ({\bf b}-{\bf Az}_{k+1} )$\\
     $k \leftarrow k + 1$
\quntil stopping criterion is satisfied.
\end{algorithm}
\vspace{0.3cm}

It is also possible (and even recommended) to update the value
of $\mu$ in each iteration \cite{NocedalWright}, \cite[Chap. 9]{Bazaraa}.
However, unlike in the quadratic penalty approach, the ALM/MM does
not require $\mu$ to be taken to infinity to guarantee convergence
to the solution of the constrained problem (\ref{constrained_linear}).

Notice that (after a straightforward complete-the-squares procedure)
the terms added to $E({\bf z})$ in the definition
of the augmented Lagrangian  ${\cal L}_A ({\bf z},\blambda_k,\mu)$ in
(\ref{augmented_L}) can be written as a single quadratic term (plus a
constant independent of ${\bf z}$, thus irrelevant for the ALM/MM),
leading to the following alternative form of the algorithm (which makes
clear its equivalence with the Bregman iterative method \cite{YinOsherGoldfarbDarbon}):

\vspace{0.3cm}
\begin{algorithm}{ALM/MM (version II)}{
\label{alg:salsa2}}
Set $k=0$, choose $\mu > 0$ and  ${\bf d}_0$.\\
\qrepeat\\
     ${\bf z}_{k+1} \in \arg\min_{{\bf z}} E({\bf z}) + \frac{\mu}{2}\|{\bf Az-d}_k\|_2^2$\\
     ${\bf d}_{k+1} = {\bf d}_{k} + ({\bf b} - {\bf Az}_{k+1})$\\
     $k \leftarrow k+1$
\quntil stopping criterion is satisfied.
\end{algorithm}
\vspace{0.3cm}

It has been shown that, with adequate initializations, the ALM/MM
generates the same sequence as a {\it proximal point algorithm} applied
to the Lagrange dual of problem (\ref{constrained_linear}) \cite{Iusem}.
Moreover, the sequence $\{{\bf d}_{k}\}$ converges to a
solution of this dual problem and all cluster points of the sequence
$\{{\bf z}_{k}\}$ are solutions of the (primal) problem (\ref{constrained_linear}) \cite{Iusem}.

\SubSection{ALM/MM for Variable Splitting}
We now show how the ALM/MM can be used to address problem
(\ref{constrained_basic}), in the particular case where $g({\bf u}) = {\bf Gu}$,
{\it i.e.},
\begin{equation}\begin{array}{cl}
 {\displaystyle \min_{{\bf u}\in \mathbb{R}^n,\, {\bf v}\in\mathbb{R}^d}} & f_1({\bf u}) + f_2({\bf v})\\
 \mbox{subject to} &{\bf G\, u} = {\bf v}, \end{array}\label{constrained_basic_linear}
\end{equation}
where ${\bf G}\in \mathbb{R}^{d\times n}.$ Problem
(\ref{constrained_basic_linear}) can be written in the form
(\ref{constrained_linear}) using the following definitions:
\begin{equation}
{\bf z} =  \left[\begin{array}{c}{\bf u}\\{\bf v}\end{array}\right],\hspace{0.5cm}
{\bf b} = {\bf 0}, \hspace{0.5cm}
{\bf A} =  [-{\bf G} \;\; {\bf I}\,],
\end{equation}
and
\begin{equation}
E({\bf z})  = f_1({\bf u}) + f_2({\bf v}).
\end{equation}
With these definitions in place, Steps 3 and 4 of the ALM/MM (version II)
can be written as follows:
\begin{eqnarray}
\left({\bf u}_{k+1} , {\bf v}_{k+1} \right) & \in & \arg\min_{{\bf u},{\bf v}}\; f_{1}({\bf u})
 + f_{2}({\bf v})  + \nonumber \\
 & & \hspace{1.18cm} \frac{\mu}{2} \|{\bf G}{\bf u} - {\bf v} - {\bf d}_k\|_2^2\label{mixed}\\
{\bf d}_{k+1} & = & {\bf d}_{k} + {\bf G}{\bf u}_{k+1} - {\bf v}_{k+1}
\end{eqnarray}

The minimization problem (\ref{mixed}) is not trivial since, in general, it
involves non-separable quadratic and possibly non-smooth terms. A natural to
address (\ref{mixed}) is to use a non-linear block-Gauss-Seidel (NLBGS) technique, in which
(\ref{mixed}) is solved by alternatingly minimizing it with respect to ${\bf u}$
and ${\bf v}$, while keeping the other variable fixed. Of course this raises
several questions: for a given ${\bf d}_k$, how much computational effort should
be spent in approximating the solution of (\ref{mixed})? Does this NLBGS
procedure converge? Experimental evidence in  \cite{GoldsteinOsher} suggests
that an efficient algorithm is obtained  by running just one NLBGS step.
It turns out that the resulting algorithm is the so-called {\it alternating
direction method of multipliers} (ADMM) \cite{EcksteinBertsekas}, which works as follows:

\vspace{0.3cm}
\begin{algorithm}{ADMM}{
\label{alg:salsa3}}
Set $k=0$, choose $\mu > 0$, ${\bf v}_0$, and  ${\bf d}_0$.\\
\qrepeat\\
   $  {\bf u}_{k+1}  \in  \arg\min_{{\bf u}} f_{1}({\bf u})
 + \frac{\mu}{2} \|{\bf G}{\bf u} - {\bf v}_k - {\bf d}_k\|_2^2$\\
  $  {\bf v}_{k+1}  \in  \arg\min_{{\bf v}} f_{2}({\bf v})
 + \frac{\mu}{2} \|{\bf G}{\bf u}_{k+1} - {\bf v} - {\bf d}_k\|_2^2$\\
     ${\bf d}_{k+1} = {\bf d}_{k} + {\bf Gu}_{k+1} - {\bf v}_{k+1}$\\
     $k \leftarrow k+1$
\quntil stopping criterion is satisfied.
\end{algorithm}
\vspace{0.3cm}

For later reference, we now recall the theorem by Eckstein and Bertsekas, in which
convergence of (a generalized version of) ADMM is shown. This theorem applies to
problems of the form (\ref{unconstrained_basic}) with  $g({\bf u}) = {\bf G}{\bf u}$,
{\it i.e.},
\begin{equation}
\min_{{\bf u}\in \mathbb{R}^n} f_1({\bf u}) + f_2\left({\bf G\, u}\right),\label{unconstrained_basic_linear}
\end{equation}
of which (\ref{constrained_basic_linear}) is the constrained optimization reformulation.

\vspace{0.3cm}
\begin{theorem}[Eckstein-Bertsekas, \cite{EcksteinBertsekas}]
\label{th:Eckstein}{\sl
Consider problem  (\ref{unconstrained_basic_linear}), where
$f_1$ and $f_2$ are closed, proper convex functions,
and ${\bf G}\in\mathbb{R}^{d\times n}$ has full column rank. Consider arbitrary
 $\mu>0$ and ${\bf v}_0, {\bf d}_0\in \mathbb{R}^d$.
 Let $\{\eta_k \geq 0, \; k=0,1,...\}$ and $\{\nu_k \geq 0,
\; k=0,1,...\}$ be two sequences such that
\[
\sum_{k=0}^\infty \eta_k < \infty \;\;\;\mbox{and} \;\;\; \sum_{k=0}^\infty \nu_k < \infty.
\]
Consider three sequences $\{{\bf u}_k \in \mathbb{R}^{n}, \; k=0,1,...\}$,
$\{{\bf v}_k \in \mathbb{R}^{d}, \; k=0,1,...\}$, and $\{{\bf d}_k \in \mathbb{R}^{d}, \; k=0,1,...\}$
that satisfy
\begin{eqnarray}
 \eta_k & \geq & \left\| {\bf u}_{k+1} - \arg\min_{{\bf u}} f_{1}({\bf u})
 + \frac{\mu}{2} \|{\bf G}{\bf u} \! - \!{\bf v}_k \! -\! {\bf d}_k\|_2^2 \right\| \nonumber\\
 \nu_k & \geq & \left\| {\bf v}_{k+1}  - \arg\min_{{\bf v}} f_{2}({\bf v})
 + \frac{\mu}{2} \|{\bf G}{\bf u}_{k+1} \! - \! {\bf v} \! - \! {\bf d}_k\|_2^2 \right\| \nonumber\\
 {\bf d}_{k+1} & = & {\bf d}_{k} + {\bf Gu}_{k+1} - {\bf v}_{k+1}.\nonumber
 \end{eqnarray}
Then, if (\ref{unconstrained_basic_linear}) has a solution, the sequence $\{{\bf u}_k\}$ converges,
${\bf u}_k \rightarrow {\bf u}^*$, where ${\bf u}^*$ is a solution of (\ref{unconstrained_basic_linear}).
If (\ref{unconstrained_basic_linear}) does not have a solution, then at least one of the
sequences $\{{\bf v}_k \}$ or $\{{\bf d}_k\}$ diverges.}
\end{theorem}
\vspace{0.3cm}

Notice that the ADMM algorithm defined above generates sequences $\{{\bf u}_k\}$, $\{{\bf v}_k\}$,
and $\{{\bf d}_k \}$ which satisfy the conditions in Theorem \ref{th:Eckstein} in a strict
sense ({\it i.e.}, with $\eta_k = \mu_k = 0$). One of the important consequences of this
theorem is that it shows that it is not necessary to exactly solve the minimizations
in lines 3 and 4 of ADMM; as long as sequence of errors is absolutely summable, convergence
is not compromised.

The proof of Theorem \ref{th:Eckstein} is based on the
equivalence between ADMM and the so-called Douglas-Rachford splitting method (DRSM)
applied to the dual of  problem (\ref{unconstrained_basic_linear}).
The DRSM was recently used for image recovery problems in \cite{CombettesPesquet}.
For recent and comprehensive reviews of ALM/MM, ADMM, DRSM, and their relationship
with Bregman and split-Bregman methods, see \cite{Esser}, \cite{Setzer}.

\Section{Proposed Method}
\label{sec:salsa}
\SubSection{Constrained Optimization Formulation of Image Recovery}
We now return to the unconstrained optimization formulation of
regularized image recovery, as defined in (\ref{problem_common}).
This problem can be written in the form (\ref{unconstrained_basic_linear}),
with
\begin{eqnarray}
f_1({\bf x}) & = & \frac{1}{2}\|{\bf Ax-y}\|_2^2 \label{eq:def_f1}\\
f_2({\bf x}) & = & \tau \phi({\bf x})\label{eq:def_f2}\\
{\bf G} & = & {\bf I}\label{eq:def_G}.
\end{eqnarray}
The constrained optimization formulation is thus
\begin{equation}\begin{array}{cl}
 {\displaystyle \min_{{\bf x},{\bf v}\in \mathbb{R}^n}} &  \frac{1}{2}\|{\bf Ax-y}\|_2^2 + \tau \phi({\bf v})\\
 \mbox{subject to} &{\bf x} = {\bf v}.
 \end{array}\label{constrained_image_linear}
\end{equation}

At this point, we are in a position to clearly explain the
difference between this formulation and the splitting exploited in
split-Bregman methods (SBM) for image recovery \cite{GoldsteinOsher}.
In those methods, the focus of attention is a non-separable
regularizer that can be written as $\phi({\bf x}) = \varphi ({\bf D\, x})$,
as is the case of the TV norm. The variable splitting used in SBM
addresses this non-separability by defining the following constrained
optimization formulation:
\begin{equation}\begin{array}{cl}
{\displaystyle \min_{{\bf x},{\bf v}\in \mathbb{R}^n}} & \frac{1}{2}\|{\bf Ax-y}\|_2^2 + \tau \varphi({\bf v})\\
\mbox{subject to} &{\bf D\, x} = {\bf v}.
\end{array}\label{constrained_image_linear_SBM}
\end{equation}

In contrast, we assume that the Moreau proximal mapping
associated to the regularizer $\phi$, {\it i.e.}, the function
$\bPsi_{\tau\phi}(\cdot)$ defined in (\ref{MPM}), can be
computed efficiently. The goal of our splitting is not
to address the difficulty raised by a non-separable and
non-quadratic regularizer, but to exploit second order (Hessian)
information of the function $f_1$, as will be shown below.

\SubSection{Algorithm and Its Convergence}
Inserting the definitions given in (\ref{eq:def_f1})--(\ref{eq:def_G}) in the
ADMM presented in the previous section yields the proposed SALSA
({\it split augmented Lagrangian shrinkage algorithm}).

\vspace{0.3cm}
\begin{algorithm}{SALSA}{
\label{alg:salsa4}}
Set $k=0$, choose $\mu > 0$, ${\bf v}_0$, and  ${\bf d}_0$.\\
\qrepeat\\
   ${\bf x}'_{k} =  {\bf v}_{k} + {\bf d}_k$ \\
   ${\bf x}_{k+1} =  \arg\min_{{\bf x}}  \|{\bf A\, x} - {\bf y}\|_2^2 + \mu \|{\bf x - x }'_k\|_2^2$\\
   ${\bf v}'_k =  {\bf x}_{k+1} - {\bf d}_k$\\
   $   {\bf v}_{k+1} =  \arg\min_{{\bf v}} \tau\phi({\bf v}) + \frac{\mu}{2}\|{\bf v}-{\bf v}'_{k}\|_2^2$\\
     $  {\bf d}_{k+1} = {\bf d}_{k} + {\bf x}_{k+1} - {\bf v}_{k+1} $\\
     $  k \leftarrow k+1$
\quntil stopping criterion is satisfied.
\end{algorithm}
\vspace{0.3cm}

Notice that SALSA is an instance of ADMM with ${\bf G=I}$; thus, the full
column rank condition on ${\bf G}$ in Theorem 1 is satisfied. If the
minimizations in lines 4 and 6 are solved exactly, we can then invoke
Theorem 1 to guarantee to convergence of SALSA.

In line 4 of SALSA, a strictly convex quadratic function has to be minimized;
which leads to the following linear system
\begin{equation}
{\bf x}_{k+1} = \left({\bf A}^H{\bf A} + \mu\, {\bf I}\right)^{-1}\left( {\bf A}^H {\bf y} +
\mu\, {\bf x}'_{k} \right).\label{eq:Wiener}
\end{equation}
As shown in the next subsection, this linear system can be solved
exactly (naturally, up to numerical precision), {\it i.e.}, non-iteratively,
for a comprehensive set of situations of interest.
The matrix ${\bf A}^H{\bf A} + \mu\, {\bf I}$ can be seen as
a regularized (by the addition of $\mu{\bf I}$) version of the
Hessian of $f_1({\bf x}) = (1/2)\|{\bf A x - y}\|_2^2$, thus SALSA
does use second order information of this function.
Notice also that (\ref{eq:Wiener}) is formally similar to the
{\it maximum a posteriori} (MAP) estimate of ${\bf x}$, from
observations $ {\bf y} = {\bf A} {\bf x} + {\bf n}$ (where ${\bf n}$
is white Gaussian noise of variance $1/\mu$) under a Gaussian prior
of mean ${\bf x}'_k$ and covariance ${\bf I}$.

The problem in line 6 is, by definition, the Moreau proximal
mapping of $\phi$ applied to ${\bf v}'_k$, thus its solution can
be written as
\begin{equation}
{\bf v}_{k+1} = \bPsi_{\tau\phi/\mu}({\bf v}'_k).
\end{equation}
If this mapping can be computed exactly in closed form, for
example, if $\phi({\bf x}) = \|{\bf x}\|_1$ thus $\bPsi$ is
simply a soft threshold, then, by Theorem 1, SALSA is guaranteed
to converge. If $\bPsi$ does not have a closed form solution and
requires itself an iterative algorithm ({\it e.g.}, if $\phi$ is
the TV norm), then convergence of SALSA still holds if one can guarantee
that the error sequence $\nu_k$ (see Theorem 1) is summable. This can
be achieved (at least approximately) if the iterative algorithm used
to approximate $\bPsi$ is initialized with the result of the previous
outer iteration, and a decreasing stopping threshold is used.

\SubSection{Computing ${\bf x}_{k+1}$}\label{sec:computingxk}
As stated above, we are interested in problems where it is not feasible to
explicitly form matrix ${\bf A}$; this might suggest that it is not easy,
or even feasible, to compute the inverse in (\ref{eq:Wiener}).
However, as shown next, in a number of problems of interest, this
inverse can be computed very efficiently.

\vspace{0.3cm}
\subsubsection{Deconvolution with Analysis Prior} In this case we have
${\bf A = B}$ (see (\ref{problem}), (\ref{problem2}), and (\ref{problem_common})),
where ${\bf B}$ is the matrix representation of a convolution.
This is the simplest case, since the inverse
$\left({\bf B}^H{\bf B} + \mu\, {\bf I}\right)^{-1}$ can be computed in the
Fourier domain. Although this is an elementary and well-known fact,
we include the derivation for the sake of completeness.
Assuming that the convolution is periodic (other boundary
conditions can be addressed with minor changes), ${\bf B}$ is a
block-circulant matrix with circulant blocks which can be factorized as
\begin{equation}
{\bf B = U}^H {\bf D U},\label{eq:diagonalization}
\end{equation}
where ${\bf U}$ is the matrix that represents the
2D discrete Fourier transform (DFT),
${\bf U}^H = {\bf U}^{-1} $ is its inverse (${\bf U}$ is unitary,
{\it i.e.}, ${\bf UU}^H = {\bf U}^H {\bf U} = {\bf I}$),
and ${\bf D}$ is a diagonal matrix containing the DFT
coefficients of the convolution operator represented by ${\bf B}$.
Thus,
\begin{eqnarray}
\left({\bf B}^H{\bf B} + \mu\, {\bf I}\right)^{-1} & = & \left({\bf U}^H{\bf D}^*{\bf D}{\bf U} + \mu\, {\bf U}^H{\bf U}\right)^{-1}\\
& = & {\bf U}^H\left(|{\bf D}|^2 + \mu\, {\bf I}\right)^{-1}{\bf U},\label{eq:Wiener2}
\end{eqnarray}
where $(\cdot)^*$ denotes complex conjugate and $|{\bf D}|^2$ the
squared absolute values of the entries of the diagonal matrix ${\bf D}$.
Since $|{\bf D}|^2 + \mu\, {\bf I}$ is diagonal, its inversion has
linear cost $O(n)$. The products by ${\bf U}$  and ${\bf U}^H$ can
be carried out with $O(n\log n)$ cost using the FFT algorithm.
The expression in (\ref{eq:Wiener2}) is a Wiener filter in the
frequency domain.

\vspace{0.3cm}
\subsubsection{Deconvolution with Frame-Based Synthesis Prior}
In this case, we have a problem of the form (\ref{problem}), {\it i.e.},
${\bf A = BW}$, thus the inversion that needs to be performed is
$\left({\bf W}^H{\bf B}^H{\bf BW} + \mu\, {\bf I}\right)^{-1}$.
Assuming that ${\bf B}$ represents a (periodic) convolution, this
inversion may be sidestepped under the assumption that matrix
${\bf W}$ corresponds to a normalized tight frame (a Parseval frame),
{\it i.e.}, ${\bf W\,W}^H = {\bf I}$. Applying the
Sherman–-Morrison-–Woodbury (SMW) matrix inversion formula   yields
{\small
\[
\left({\bf W}^H{\bf B}^H{\bf BW} + \mu\, {\bf I}\right)^{-1} \!\! =
\frac{1}{\mu} \, ( {\bf I} -   {\bf W}^H \underbrace{{\bf B}^H \left({\bf B} {\bf B}^H \!
+ \mu\, {\bf I}\right)^{-1} \!\! {\bf B}}_{\bf F} {\bf W}).
\] } Let's focus on the term ${\bf F \equiv  B}^H \left({\bf B} {\bf B}^H \!
+ \mu\, {\bf I}\right)^{-1}{\bf B}$; using the factorization
(\ref{eq:diagonalization}), we have
\begin{equation}
{\bf F} =
 {\bf U}^H {\bf D^*} \left( |{\bf D}|^2  + \mu\, {\bf I}\right)^{-1} {\bf D U}.
 \label{eq:fft_filter}
\end{equation}
Since all the matrices in
${\bf D^*} \left( |{\bf D}|^2  + \mu\, {\bf I}\right)^{-1} {\bf D}$ are diagonal,
this expression can be computed with $O(n)$ cost, while the products by
${\bf U}$ and ${\bf U}^H$ can be computed with
$O(n\log n)$ cost using the FFT. Consequently, products by matrix ${\bf F}$
(defined in (\ref{eq:fft_filter})) have $O(n\log n)$ cost.

Defining ${\bf r}_k = \left( {\bf A}^H {\bf y} + \mu\, {\bf x}'_{k} \right)=
\left( {\bf W}^H {\bf B}^H {\bf y} + \mu\, {\bf x}'_{k} \right)$, allows writing
(\ref{eq:Wiener}) compactly as
\begin{equation}
{\bf x}_{k+1} = \frac{1}{\mu} \, \left({\bf r}_k - {\bf W}^T {\bf F} {\bf W}\; {\bf r}_k\right) .
\label{GS_linear_5}
\end{equation}
Notice that multiplication by ${\bf F}$ corresponds to applying an image filter in
the Fourier domain. Finally, notice also that the term ${\bf B}^H{\bf W}^H {\bf y}$
can be precomputed, as it doesn't change during the algorithm.

The leading cost of each application of (\ref{GS_linear_5}) will be either $O(n\log n)$
or the cost of the products by ${\bf W}^H$ and ${\bf W}$. For most tight frames
used in image processing, these products correspond to direct and inverse
transforms for which fast algorithms exist. For example, when ${\bf W}^H$ and
${\bf W}$ are the inverse and direct translation-invariant wavelet transforms,
these products can be computed using the undecimated wavelet transform with
$O(n\log n)$ total  cost \cite{lang}. Curvelets also constitute a Parseval frame
for which fast $O(n\log n)$ implementations of the forward and inverse transform exist
\cite{CandesDemanetDonohoYing}. Yet another example of a redundant Parseval
frame is the complex wavelet transform, which has $O(n)$ computational cost
\cite{Kingsbury}, \cite{Selesnick}. In conclusion, for a large class of choices
of ${\bf W}$, each iteration of the SALSA algorithm has $O(n\log n)$ cost.

\vspace{0.3cm}
\subsubsection{Missing Pixels: Image Inpainting}
In the analysis prior case (TV-based), we have ${\bf A} = {\bf B}$,
where the observation matrix ${\bf B}$ models the loss of some
image pixels. Matrix ${\bf A=B}$ is thus an $m\times n$ binary
matrix, with $m<n$, which can be obtained by taking a subset of
rows of an identity matrix. Due to its particular structure, this matrix
satisfies ${\bf B}{\bf B}^T = {\bf I}$. Using this fact together with
the SMW formula leads to
\begin{eqnarray}
\left({\bf B}^T{\bf B} + \mu {\bf I}\right)^{-1} = \frac{1}{\mu}\left(
{\bf I} - \frac{1}{1+\mu}\, {\bf B}^T{\bf B}\right).\label{eq:mask_matrix}
\end{eqnarray}
Since ${\bf B}^T{\bf B}$ is equal to an identity matrix with some
zeros in the diagonal (corresponding to the positions of the
missing observations), the matrix in (\ref{eq:mask_matrix}) is
diagonal with elements either equal to $1/(\mu +1)$ or $1/\mu$.
Consequently,  (\ref{eq:Wiener}) corresponds simply to multiplying
$({\bf B}^H {\bf y} + \mu {\bf x}'_k)$ by this diagonal matrix, which is
an $O(n)$ operation.

In the synthesis prior case, we have ${\bf A} = {\bf BW}$, where ${\bf B}$
is the binary sub-sampling matrix defined in the previous paragraph.
Using  the SMW formula yet again, and the fact that
${\bf B}{\bf B}^T = {\bf I}$, we have
{\small
\begin{equation}
\left({\bf W}^H{\bf B}^H{\bf BW} + \mu\, {\bf I}\right)^{-1} \!\! =
\frac{1}{\mu} \, {\bf I} -  \frac{\mu}{1+\mu}\, {\bf W}^H {\bf B}^T {\bf B} {\bf W} .
\label{eq:missing_synth}
\end{equation} }
As noted in the previous paragraph, ${\bf B}^T{\bf B}$ is equal to an identity matrix with
zeros in the diagonal (corresponding to the positions of the
missing observations), {\it i.e.}, it is a binary mask. Thus, the multiplication
by  ${\bf W}^H {\bf B}^H {\bf B} {\bf W} $ corresponds to synthesizing
the image, multiplying it by this mask, and computing the representation
coefficients of the result. In conclusion, the cost of (\ref{eq:Wiener})
is again that of the products by ${\bf W}$ and ${\bf W}^H$, usually $O(n\log n)$.

\vspace{0.3cm}
\subsubsection{Partial Fourier Observations: MRI Reconstruction.}
The final case considered is that of partial Fourier observations,
which is used to model magnetic resonance image (MRI) acquisition \cite{Lustig},
and has been the focus of much recent interest due to its connection
to compressed sensing \cite{Candes,ec1,donoho1}.
In the TV-regularized case, the observation matrix has the form ${\bf A}={\bf B}{\bf U}$,
where ${\bf B}$ is an $m\times n$ binary matrix, with $m<n$, similar to the one
in the missing pixels case (it is formed by a subset of rows of an identity matrix),
and ${\bf U}$ is the DFT matrix. This case is similar to (\ref{eq:missing_synth}),
with ${\bf U}$ and ${\bf U}^H$ instead of ${\bf W}$ and ${\bf W}^H$,
respectively. The cost of (\ref{eq:Wiener}) is again that of the products by ${\bf U}$ and ${\bf U}^H$,
{\it i.e.}, $O(n\log n)$ if we use the FFT.

In the synthesis case, the observation matrix has the form ${\bf A}={\bf B}{\bf U}{\bf W}$.
Clearly, the case is again similar to (\ref{eq:missing_synth}),
but with ${\bf UW}$ and ${\bf W}^H{\bf U}^H$ instead of ${\bf W}$ and ${\bf W}^H$,
respectively. Again, the cost of (\ref{eq:Wiener}) is $O(n\log n)$, if the FFT is
used to compute the products by ${\bf U}$ and ${\bf U}^H$ and fast frame transforms
are used for the products by ${\bf W}$ and ${\bf W}^H$.

\Section{Experiments}
\label{sec:experiments}

In this section, we report results of experiments aimed at comparing
the speed of SALSA with that of the current state of the art methods
(all of which are freely available online): TwIST\footnote{Available at \url{ http://www.lx.it.pt/~bioucas/code/TwIST_v1.zip}} \cite{TwIST},
SpaRSA\footnote{Available at \url{http://www.lx.it.pt/~mtf/SpaRSA/}}
\cite{SpaRSA_SP},
and FISTA\footnote{Available at \url{http://iew3.technion.ac.il/~becka/papers/wavelet_FISTA.zip}}
\cite{FISTA}. We consider three standard and often studied imaging inverse problems:
image deconvolution (using both wavelet and TV-based regularization); image restoration from
missing samples (inpainting); image reconstruction from partial Fourier observations,
which (as mentioned above) has been the focus of much recent interest due to its connection
with compressed sensing and the fact that it models MRI acquisition \cite{Lustig}.
All experiments were performed using MATLAB for Windows XP, on a desktop computer
equipped with an Intel Pentium-IV $3.0$ GHz processor and $1.5$GB of  RAM.
To compare the speed of the algorithms, in a way that is as independent as possible
from the different stopping criteria, we first run SALSA and then the other algorithms
until they reach the same value of the objective function. The value of $\mu$ for
fastest convergence was found to differ (though not very much) in each case, but a
good rule of thumb, adopted in all the experiments, is $\mu = 0.1\tau$.

\begin{table}[hbt]
\centering
\caption{Details of the image deconvolution experiments.}\label{decon_problems}
\begin{tabular}{|c|l|l|}
  \hline
 Experiment &  blur kernel  \rule[-0.1cm]{0cm}{0.4cm}  & $\sigma^2$ \\ \hline
1  &$9\times 9$ uniform & $0.56^2$ \\
2A & Gaussian & 2\\
2B & Gaussian & 8\\
3A  & $h_{ij} = 1/(1 + i^2 + j^2)$ & 2 \\
3B  & $h_{ij} = 1/(1 + i^2 + j^2)$ & 8 \\
\hline
\end{tabular}
\end{table}

\subsection{Image Deblurring with wavelets}
We consider five benchmark deblurring problems \cite{FigueiredoNowak2003},
summarized in Table~\ref{decon_problems}, all on the well-known Cameraman image.
The regularizer is $\phi(\bbeta) = \|\bbeta\|_1$, thus $\bPsi_{\tau\phi}$ is an
element-wise soft threshold. The blur operator ${\bf B}$ is applied via the FFT.
The regularization parameter $\tau$ is hand tuned in each case for best improvement
in SNR, so that the comparison is carried out in the regime that is relevant in
practice. Since the restored images are visually indistinguishable from those
obtained in \cite{FigueiredoNowak2003}, and the SNR improvements are also very
similar, we simply report computation times.

\begin{figure}[tp]
\centering
\subfigure[]{
\includegraphics[width=0.3\textwidth]{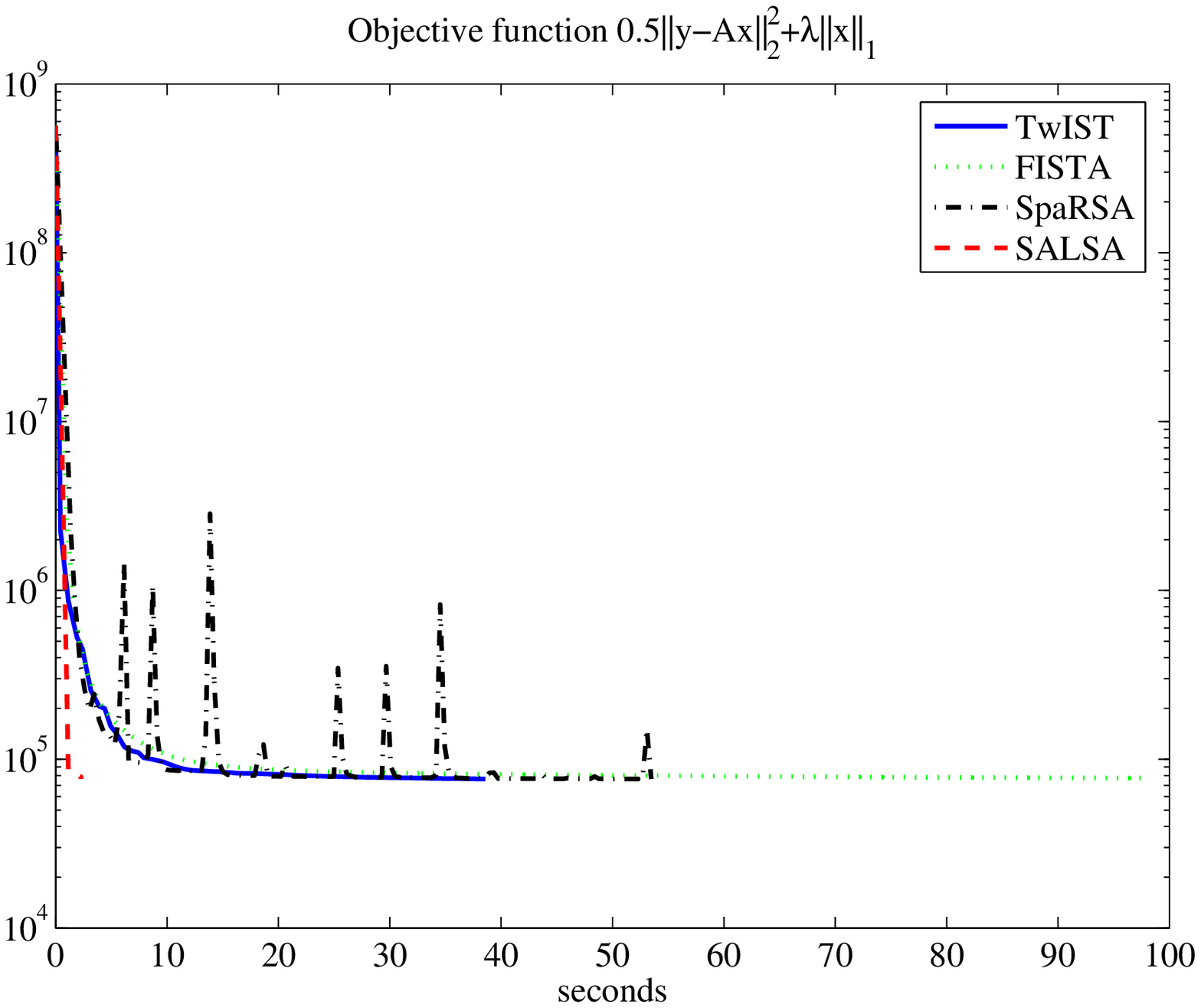}
}
\subfigure[]{
\includegraphics[width=0.3\textwidth]{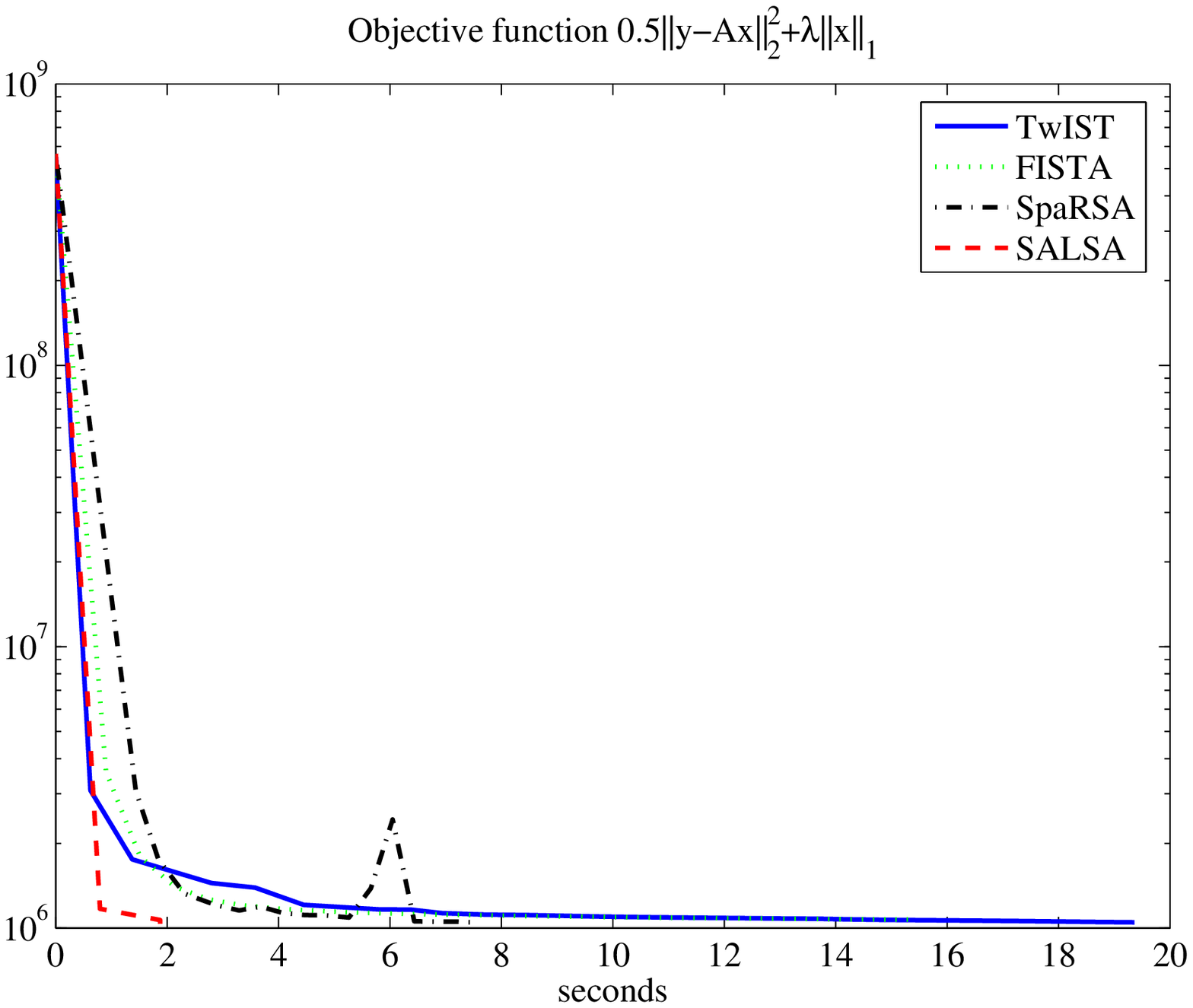}
}
\subfigure[]{
\includegraphics[width=0.3\textwidth]{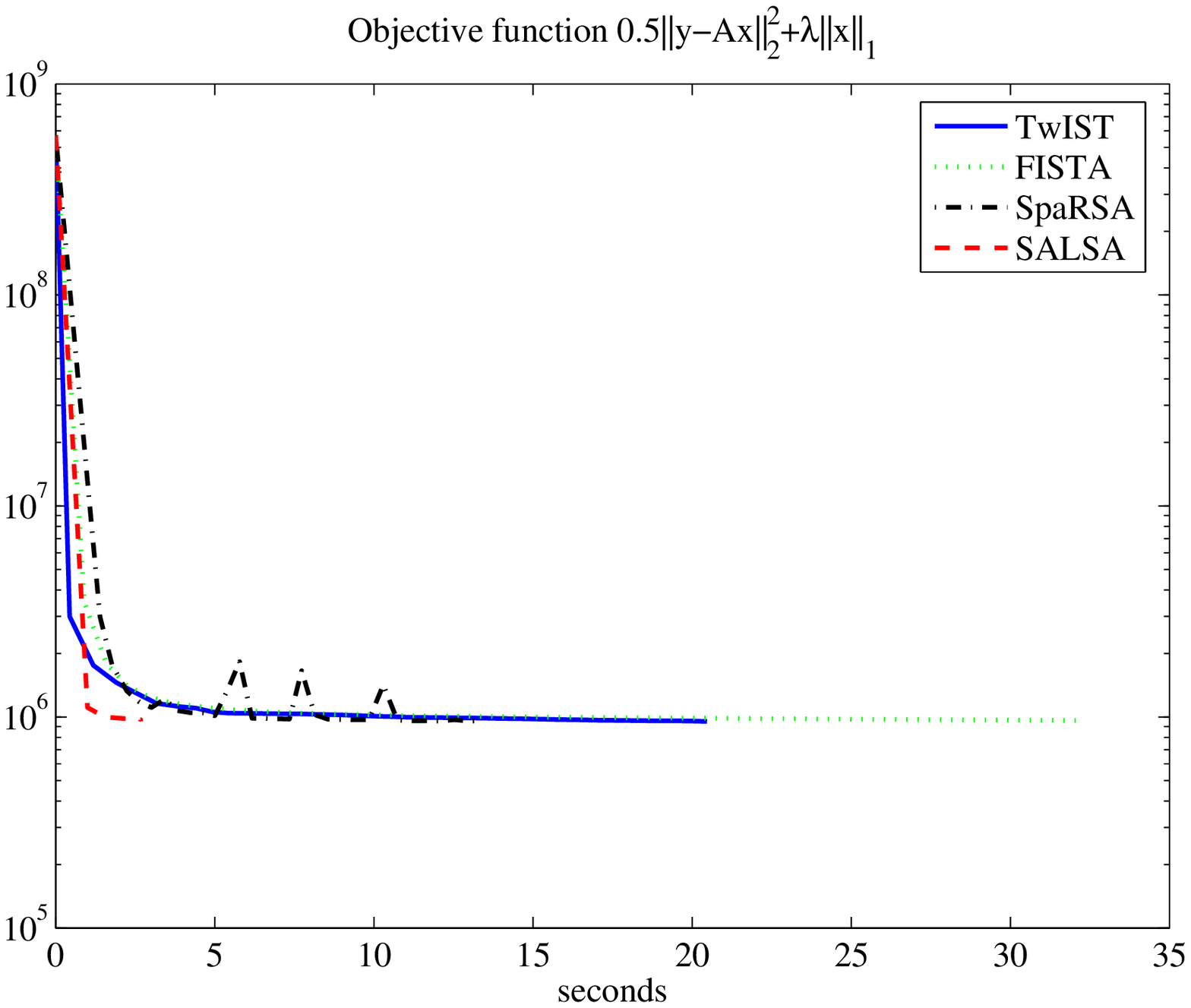}
}
\caption{Objective function evolution (redundant wavelets):  (a) experiment 1A; (b) experiment 2B; (c) experiment 3A.}
\label{fig:evolutionobjectivel1_redundant}
\end{figure}

\begin{figure}[t!]
\centering
\subfigure[]{
\includegraphics[width=0.3\textwidth]{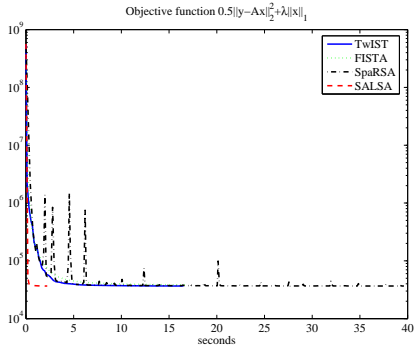}
}
\subfigure[]{
\includegraphics[width=0.3\textwidth]{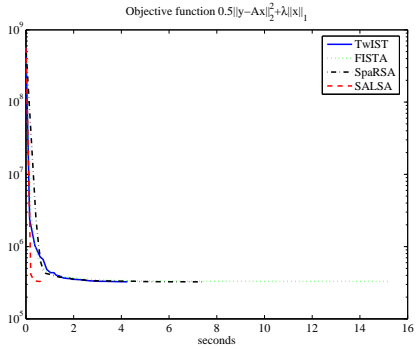}
}
\subfigure[]{
\includegraphics[width=0.3\textwidth]{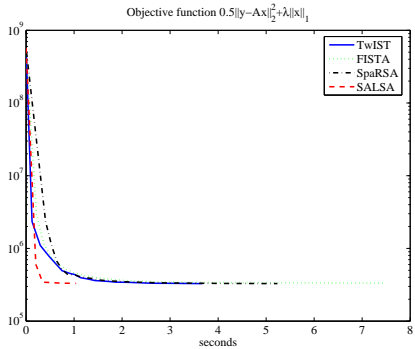}
}
\caption{Objective function evolution (orthogonal wavelets):  (a) experiment 1A; (b) experiment 2B; (c) experiment 3A.}
\label{fig:evolutionobjectivel1_orthogonal}
\end{figure}

In the first set of experiments, ${\bf W}$ is a redundant Haar wavelet frame
with four levels. The CPU times taken by each of the algorithms are presented in Table~\ref{tab:deconl1redundant_times}. In the second set of experiments, ${\bf W}$ is an orthogonal Haar wavelet basis; the results are reported in Table~\ref{tab:deconl1orthogonal_times}. To visually illustrate the relative speed of the algorithms, Figures~\ref{fig:evolutionobjectivel1_redundant}
and \ref{fig:evolutionobjectivel1_orthogonal} plot the evolution of the objective function (see Eq. (\ref{problem})),
versus time, in experiments $1$, $2$B, and $3$A,
for redundant and orthogonal wavelets, respectively.

%
%


\begin{table}[hbt]
\centering \caption{Image deblurring  with redundant wavelets: CPU times (in seconds).}
\label{tab:deconl1redundant_times}
\begin{tabular}{|c|l|l|l|l|}
  \hline
 Experiment & TwIST   \rule[-0.1cm]{0cm}{0.4cm}  & SpARSA  & FISTA & SALSA \\
\hline
1 & 38.5781 & 53.4844 & 98.2344 & 2.26563 \\
2A & 33.8125 & 42.7656 & 65.3281 & 4.60938 \\
2B & 35.2031 & 70.7031 & 112.109 & 12.0313 \\
3A & 20.4688 & 13.3594 & 32.2969 & 2.67188 \\
3B & 9.0625 & 5.8125 & 18.0469 & 2.07813 \\
\hline
\end{tabular}
\end{table}

\begin{table}[h!]
\centering \caption{Image deblurring with orthogonal wavelets: CPU times (in seconds).}
\label{tab:deconl1orthogonal_times}
\begin{tabular}{|c|l|l|l|l|}
  \hline
 Experiment & TwIST   \rule[-0.1cm]{0cm}{0.4cm}  & SpARSA  & FISTA & SALSA \\
\hline
1 & 16.5156 & 39.6094 & 16.8281 & 2.23438 \\
2A & 10.1406 & 16.3438 & 15.9531 & 1.375 \\
2B & 5.10938 & 7.96875 & 5.3125 & 0.640625 \\
3A & 3.67188 & 5.23438 & 7.46875 & 1.03125 \\
3B & 2.57813 & 2.64063 & 3.625 & 0.5625 \\
\hline
\end{tabular}
\end{table}


\subsection{Image Deblurring with Total Variation}\label{sec:TVreg_exp}
The same five image deconvolution problems listed in Table~\ref{decon_problems} were
also addressed using total variation (TV) regularization (more specifically, the
isotropic discrete total variation, as defined in \cite{Chambolle}). The corresponding
Moreau proximal mapping is computed using $5$ iterations of Chambolle's algorithm \cite{Chambolle}.

\begin{figure}[t!]
\centering
\subfigure[]{
\includegraphics[width=0.3\textwidth]{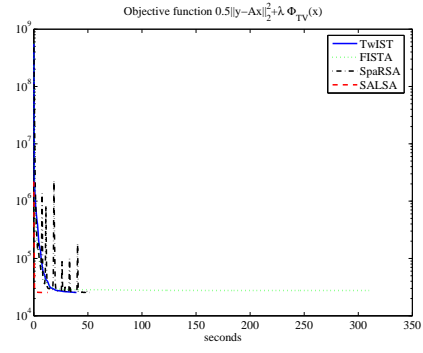}
}
\subfigure[]{
\includegraphics[width=0.3\textwidth]{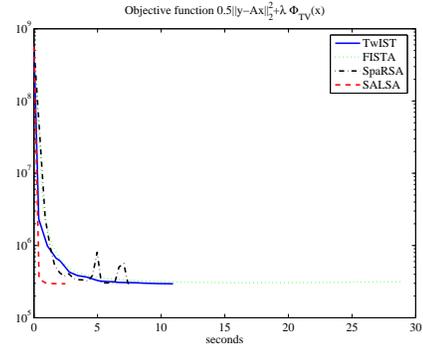}
}
\subfigure[]{
\includegraphics[width=0.3\textwidth]{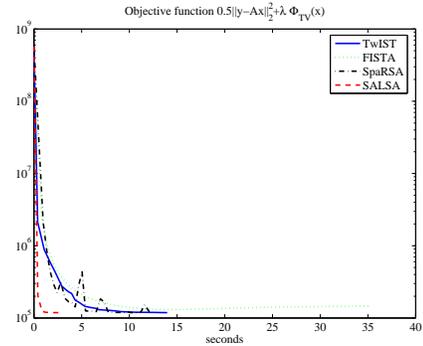}
}
\caption{Image deblurring with TV regularization - Objective function evolution: (a) $9\times 9$ uniform
blur, $\sigma=0.56$; (b) Gaussian blur, $\sigma^2=8$; (c) $h_{ij} = 1/(1 + i^2 + j^2)$ blur, $\sigma^2=2$.}
\label{fig:evolutionobjectiveTV}
\end{figure}

The CPU times taken by SALSA, TwIST, SpaRSA, and FISTA are listed in Table~\ref{tab:deconTV_times}. The evolutions of the objective functions (for experiments $1$, $2$B, and $3$A) are plotted in Figure~\ref{fig:evolutionobjectiveTV}.

We can conclude from Tables \ref{tab:deconl1redundant_times}, \ref{tab:deconl1orthogonal_times},
and \ref{tab:deconTV_times} that, in image deconvolution problems, both with wavelet-based and TV-based
regularization, SALSA is always clearly faster than the fastest of the other
competing algorithms.

\begin{table}[h!]
\centering \caption{TV-Based Image deblurring: CPU Times (in seconds).}
\label{tab:deconTV_times}
\begin{tabular}{|c|l|l|l|l|}
  \hline
 Experiment & TwIST   \rule[-0.1cm]{0cm}{0.4cm}  & SpARSA  & FISTA & SALSA \\
\hline
1 & 63.2344 & 80.0469 & 346.734 & 11.2813 \\
2A & 19.1563 & 24.1094 & 34.1406 & 4.79688 \\
2B & 10.9375 & 7.75 & 29.0156 & 2.46875 \\
3A & 13.9688 & 12.4375 & 35.2969 & 2.79688 \\
3B & 10.9531 & 7.75 & 28.3438 & 2.78125 \\
\hline
\end{tabular}
\end{table}



\subsection{MRI Image Reconstruction}
\begin{figure}[h]
\centering
\subfigure[]{
\includegraphics[width=0.25\textwidth]{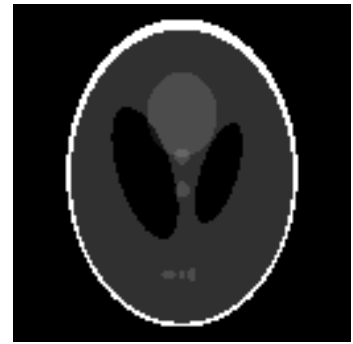}
\label{fig:phantom128}
}

\subfigure[]{
\includegraphics[width=0.25\textwidth]{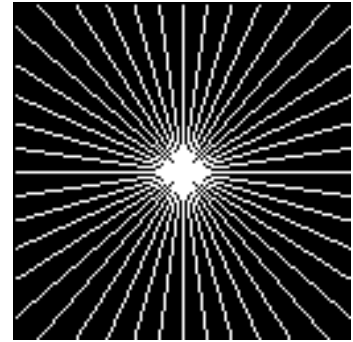}
\label{fig:MRImask22}
}

\subfigure[]{
\includegraphics[width=0.25\textwidth]{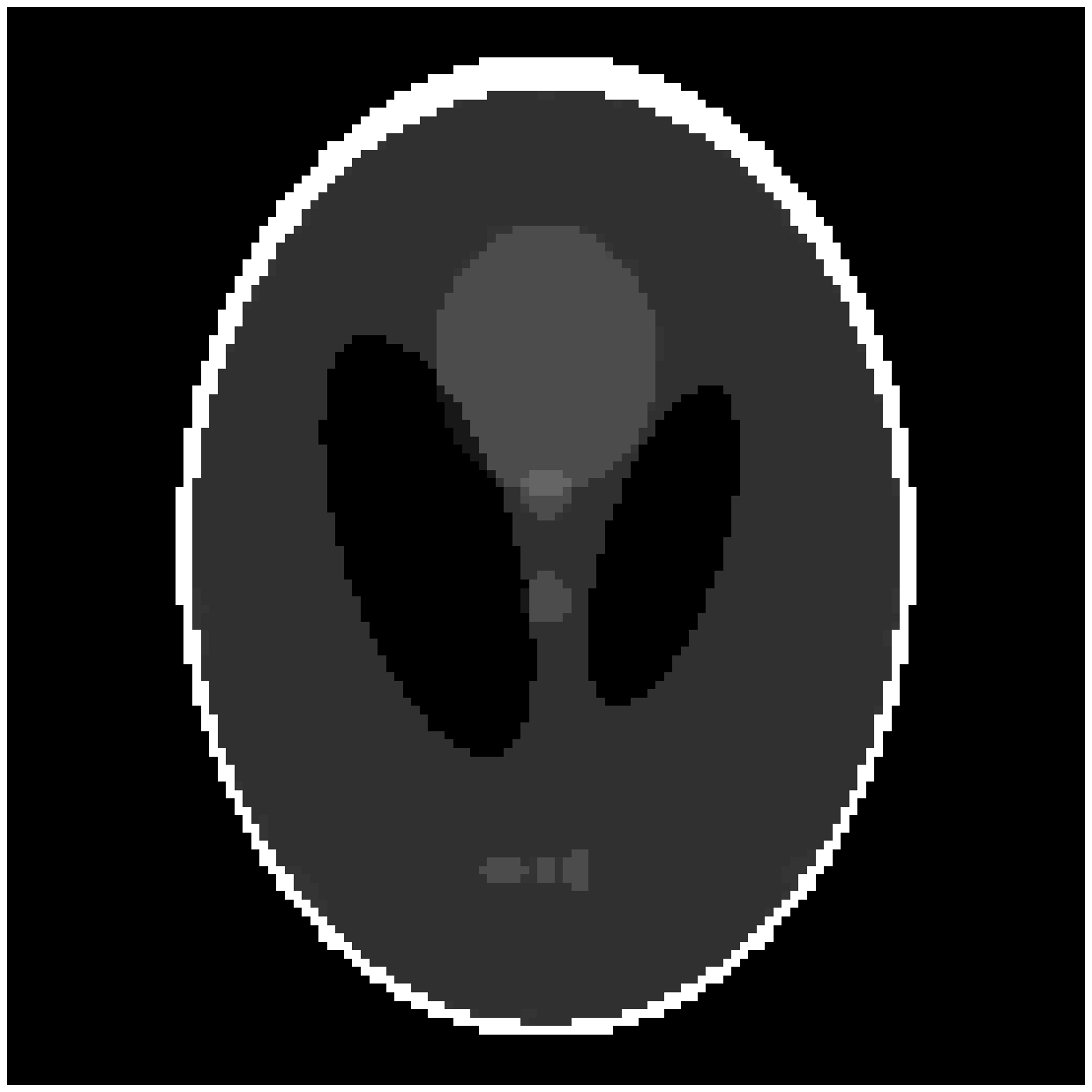}
\label{fig:estimateSALSA_mri}
}
\caption{MRI reconstruction: (a)$128\times 128$ Shepp Logan phantom; (b) Mask with 22 radial lines;  (c) image estimated using SALSA.}
\end{figure}

We consider the problem of reconstructing the $128\times 128$ Shepp-Logan phantom
(shown in Figure~\ref{fig:phantom128}) from a limited number of radial lines (22, in our
experiments, as shown in Figure~\ref{fig:MRImask22}) of its 2D discrete Fourier
transform.  The projections are also corrupted with circular complex Gaussian noise,
with variance $\sigma^2\ =\ 0.5\times 10^{-3}$. We use TV regularization (as described
in Subsection \ref{sec:TVreg_exp}), with the corresponding Moreau proximal mapping
implemented by  $40$ iterations of Chambolle's algorithm \cite{Chambolle}.


\begin{figure}[h]
\centering
\includegraphics[width=0.3\textwidth]{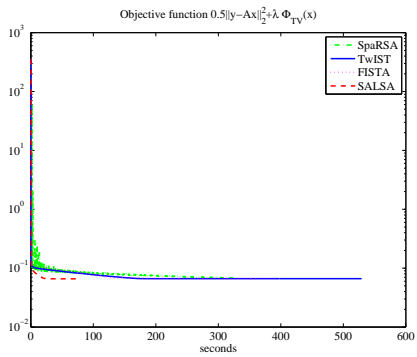}
\caption{MRI reconstruction: evolution of the objective function over time.}
\label{fig:objective_mri}
\end{figure}

Table~\ref{tab:MRI_comparison} shows the CPU times, while Figure~\ref{fig:objective_mri} plots
the evolution of the objective function over time.
Figure~\ref{fig:estimateSALSA_mri} shows the estimate obtained using SALSA (the others are, naturally,
visually indistinguishable). Again, we may conclude that SALSA is considerably faster than
the other three algorithms, while achieving comparable values of mean squared error of
the reconstructed image.

\begin{table}[hbt]
\centering \caption{MRI reconstruction: Comparison of the various algorithms.}
\label{tab:MRI_comparison}
\begin{tabular}{|c|l|l|l|l|}
  \hline
  & TwIST   \rule[-0.1cm]{0cm}{0.4cm}  & SpARSA  & FISTA & SALSA \\
\hline
Iterations & 1002 & 1001 & 1000 & 53 \\
CPU time (seconds) & 529.297 & 328.688 & 390.75 & 76.5781 \\
MSE & 4.384e-7 & 6.033e-5 & 4.644e-7 & 5.817e-7 \\
\hline
\end{tabular}
\end{table}

\subsection{Image Inpainting}
\begin{figure}[h!]
\centering
\subfigure[]{
\includegraphics[width=0.25\textwidth]{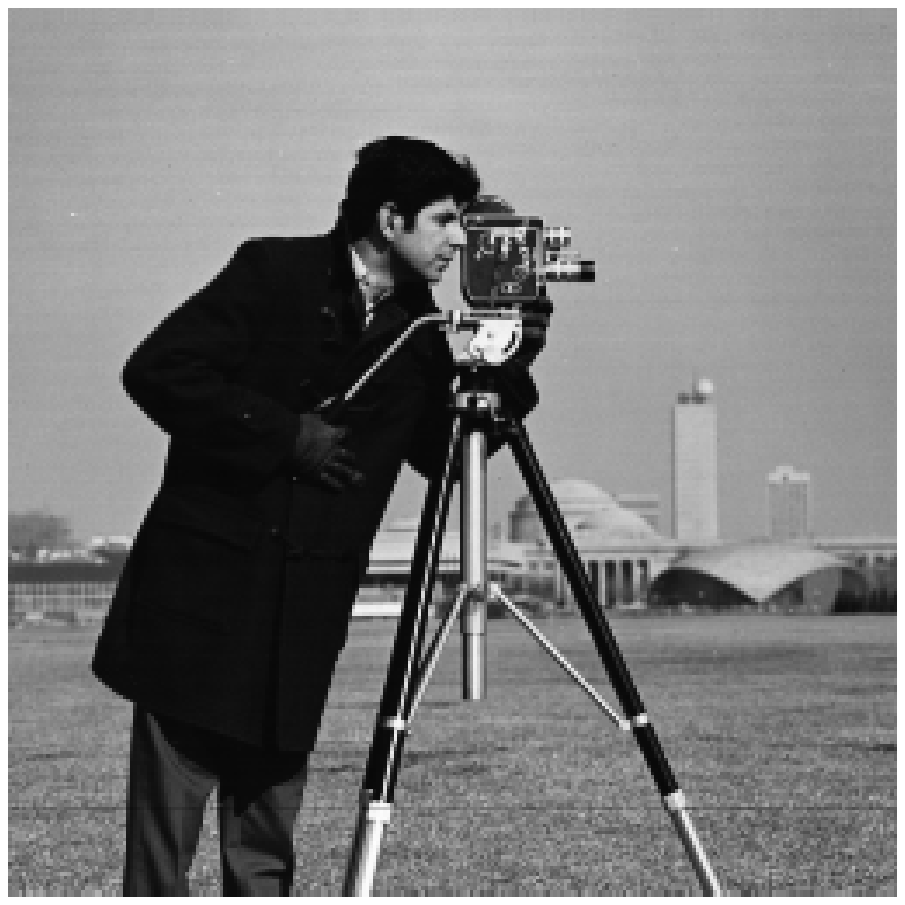}
}
\subfigure[]{
\includegraphics[width=0.25\textwidth]{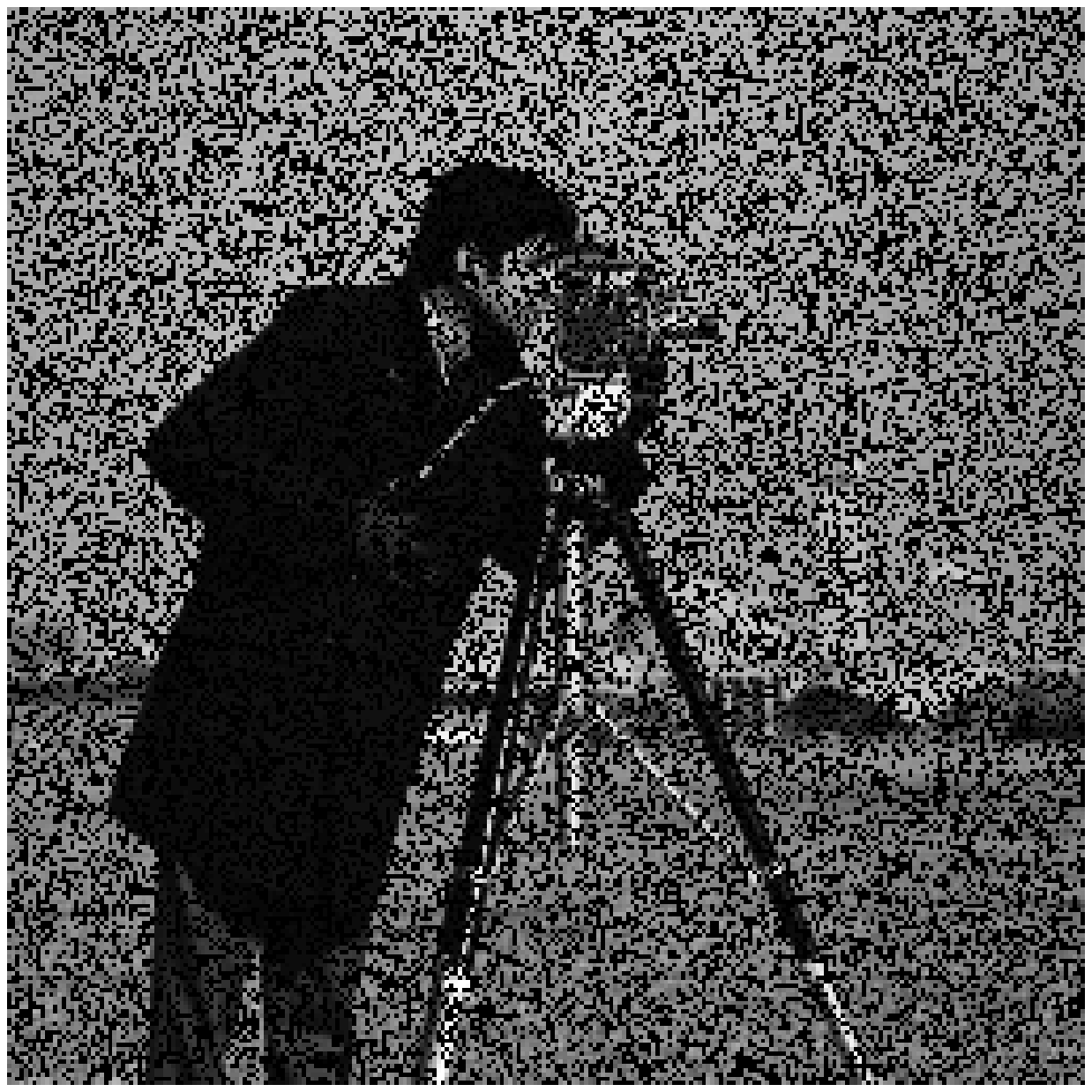}
}
\subfigure[]{
\includegraphics[width=0.25\textwidth]{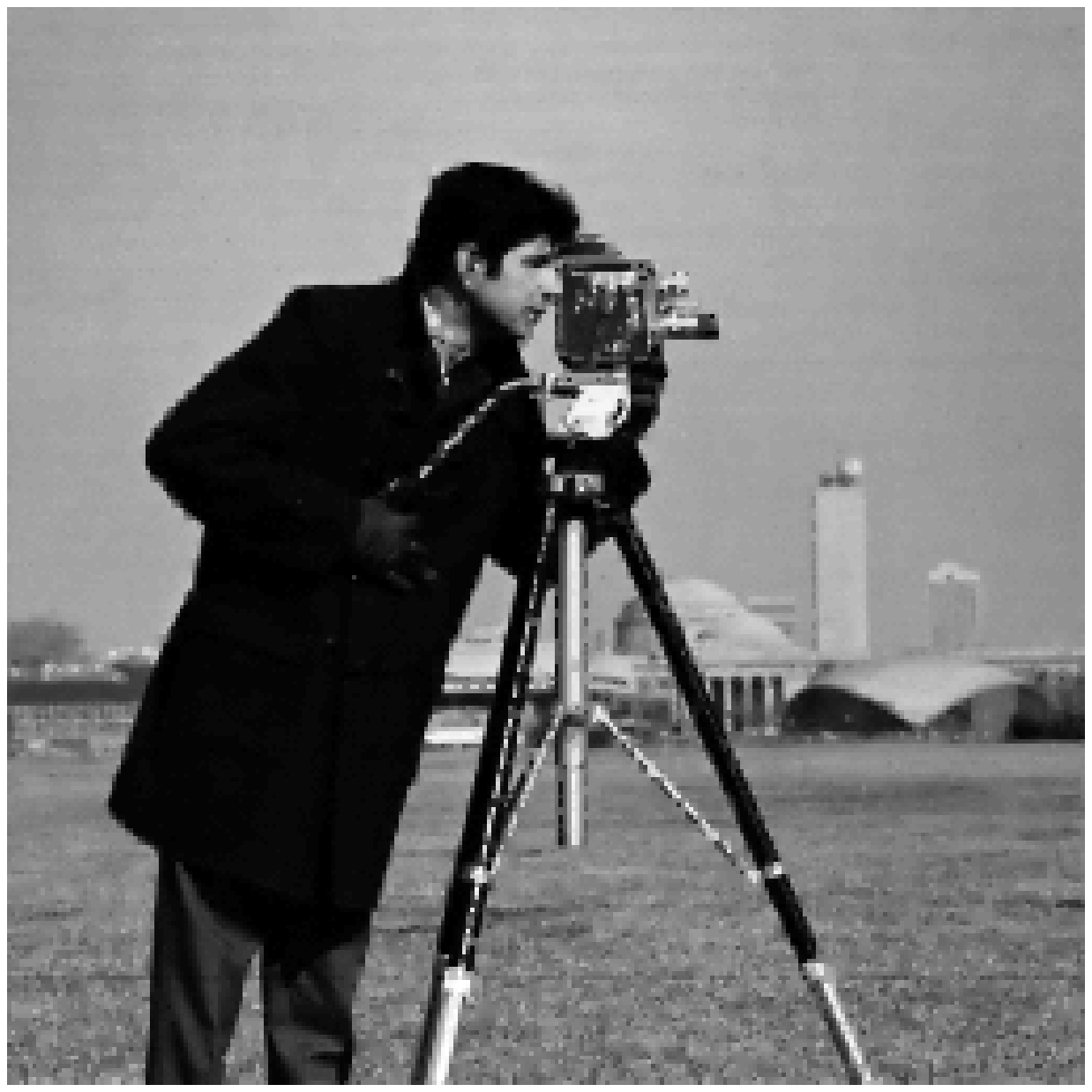}
}
\caption{Image inpainting with TV regularization: (a) Original cameraman image; (b) Image with $40\%$ pixels missing; (c) Estimated using SALSA.}
\label{fig:missingpixels}
\end{figure}

Finally, we consider an image inpainting problem, as explained in Section~\ref{sec:computingxk}.
The original image is again the Cameraman, and the observation consists in loosing $40\%$ of its pixels,
as shown in Figure~\ref{fig:missingpixels}. The observations are also corrupted with Gaussian noise
(with an SNR of $40$ dB). The regularizer is again TV implemented by $20$ iterations of Chambolle's
algorithm.

\begin{figure}[th]
\centering
\includegraphics[width=0.3\textwidth]{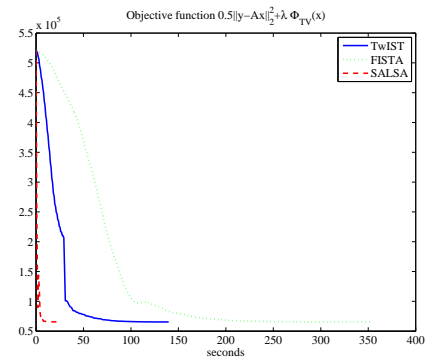}
\caption{Image inpainting: evolution of the objective function over time.}
\label{fig:objective_missingpixels}
\end{figure}

The image estimate obtained by SALSA is shown in Figure~\ref{fig:missingpixels}, with the original
also shown for comparison. The estimates obtained using TwIST and FISTA were visually very similar. Table~\ref{tab:missingdata_comparison} compares the performance of SALSA with that of TwIST and FISTA and Figure~\ref{fig:objective_missingpixels} shows the evolution of the objective function for each of the
algorithms. Again, SALSA is considerably faster than the alternative algorithms.

\begin{table}[hbt]
\centering \caption{Image inpainting: Comparison of the various algorithms.}
\label{tab:missingdata_comparison}
\begin{tabular}{|c|l|l|l|}
  \hline
  & TwIST   \rule[-0.1cm]{0cm}{0.4cm}  & FISTA & SALSA \\
\hline
Iterations & 302 & 300 & 33 \\
CPU time (seconds) & 305 & 228 & 23 \\
MSE & 105 & 101 & 99.1 \\
ISNR (dB) & 18.4 & 18.5 & 18.6 \\
\hline
\end{tabular}
\end{table}

\Section{Conclusions}
\label{sec:conclusions}

We have presented a new algorithm for solving the unconstrained
optimization formulation of regularized image reconstruction/restoration.
The approach, which can be used with different types of regularization
(wavelet-based, total variation), is based on a variable splitting technique
which yields an equivalent constrained problem. This constrained
problem is then addressed using an augmented Lagrangian method,
more specifically, the alternating direction method of multipliers (ADMM).
The algorithm uses a regularized version of the Hessian of the $\ell_2$
data-fidelity term, which can be computed efficiently for several classes
of problems. Experiments on a set of standard image recovery problems
(deconvolution, MRI reconstruction, inpainting) have shown that the
proposed algorithm (termed SALSA, for {\it split augmented Lagrangian
shrinkage algorithm}) is faster than previous state-of-the-art methods.
Current and future work involves using a similar approach to solve
constrained formulations of the forms
(\ref{eq:bp}) and (\ref{eq:bp2}).

\end{document}